\documentclass[11pt]{article}
\usepackage{amsthm,amsmath,latexsym,amssymb}

\newcommand{\bP}{{\rm |\kern-.15em P}}
\newcommand{\Q}{\kern.3em\rule{.07em}{.65em}\kern-.3em{\rm Q}}
\newcommand{\R}{{\rm I\kern-.15em R}}
\newcommand{\D}{{\rm |\kern-.15em D}}
\newcommand{\h}{{\rm |\kern-.15em H}}
\newcommand{\C}{\kern.3em\rule{.07em}{.65em}\kern-.3em{\rm C}}
\newcommand{\T}{{\rm T\kern-.35em T}}

\theoremstyle{plain}
\newtheorem{theorem}{Theorem}[section]

\theoremstyle{definition}
\newtheorem{definition}[theorem]{Definition}

\theoremstyle{remark}

\begin{document}
\title{The Jacobian Conjecture for the space of all the inner functions}
\author{Ronen Peretz}
 
\maketitle

\begin{abstract}
We prove the Jacobian Conjecture for the space of all the inner functions
in the unit disc. 
\end{abstract}

\section{Known facts}

\begin{definition}
Let $B_F$ be the set of all the finite Blaschke products defined on the unit disc 
$\mathbb{D}=\{z\in\mathbb{C}\,|\,|z|<1\}$.
\end{definition}
\noindent
{\bf Theorem A.} $f(z)\in B_F\Leftrightarrow\,\exists\,n\in\mathbb{Z}^+\,\,{\rm such}\,\,{\rm that}\,\,\forall\,
w\in\mathbb{D}$ the equation $f(z)=w$ has exactly $n$ solutions $z_1,\ldots,z_n$ in $\mathbb{D}$, counting
multiplicities. \\
\\
That follows from \cite{fm} on the bottom of page 1. \\
\\
{\bf Theorem B.} $(B_F,\circ)$ is a semigroup under composition of mappings. \\
\\
That follows by Theorem 1.7 on page 5 of \cite{kr}. \\
\\
{\bf Theorem C.} If $f(z)\in B_F$ and if $f'(z)\ne 0\,\,\forall\,z\in\mathbb{D}$ then
$$
f(z)=\lambda\frac{z-\alpha}{1-\overline{\alpha}z}
$$
for some $\alpha\in\mathbb{D}$ and some unimodular $\lambda$, $|\lambda|=1$, i.e. $f\in {\rm Aut}(\mathbb{D})$. \\
\\
For that we can look at Remark 1.2(b) on page 2, and remark 3.2 on page 14 of \cite{kr}. Also we can
look at Theorem A on page 3 of \cite{kr1}.

\section{introduction}

We remark that the last theorem (Theorem C) could be thought of, as the (validity of) Jacobian Conjecture for $B_F$. This result
is, perhaps, not surprising in view of the characterization in Theorem A above of members of $B_F$ (This is, in fact
Theorem B on page 2 of \cite{fm}. This result is due to Fatou and to Rado). For in the classical Jacobian
Conjecture one knows of a parallel result, namely: \\
\\
If $F\in{\rm et}(\mathbb{C}^2)$ and if $d_F(w)=|\{ z\in\mathbb{C}^2\,|\,F(z)=w\}|$ is a constant $N$ (independent
of $w\in\mathbb{C}^2$), then $F\in{\rm Aut}(\mathbb{C}^2)$ (because $F$ is a proper mapping). \\
Thus we are led to the following,
\begin{definition}
Let $V_F$ be the set of all holomorphic $f\,:\,\mathbb{D}\rightarrow\mathbb{D}$, such that $\exists\,N=N_f\in\mathbb{Z}^+$
(depending on $f$) for which $d_f(w)=|\{z\in\mathbb{D}\,|\,f(z)=w\}|$, $w\in\mathbb{D}$, satisfies $d_f(w)\le N_f$
$\forall\,w\in\mathbb{D}$.
\end{definition}
\noindent
We ask if the following is true:
$$
f\in V_F,\,\,f'(z)\ne 0\,\forall\,z\in\mathbb{D}\Rightarrow f\in{\rm Aut}(\mathbb{D}).
$$
The answer is negative. For example, we can take $f(z)=z/2$. So we modify the question:
$$
f\in V_F,\,\,f'(z)\ne 0\,\forall\,z\in\mathbb{D}\Rightarrow f(z)\,\,{\rm is}\,\,{\rm injective}.
$$
This could be written, alternatively as follows:
$$
f\in V_F,\,\,f'(z)\ne 0\,\forall\,z\in\mathbb{D}\Rightarrow \forall\,w\in\mathbb{D},\,\,d_f(w)\le 1.
$$
Also the answer to this question is negative. For we can take $f(z)=10^{-10}e^{10z}$ which will satisfy
the condition $f(\mathbb{D})\subset\mathbb{D}$ because of the tiny factor $10^{-10}$, while clearly
$f\in V_F$ and $f'(z)\ne 0$ $\forall\,z\in\mathbb{D}$. But $d_f(w)$ can be as large as
$$
\left[\frac{2}{2\pi/10}\right]=\left[\frac{10}{\pi}\right]=3.
$$
Thus we again need to modify the question (in order to get a more interesting result). It is not clear if the right assumption 
should include surjectivity or almost surjectivity. Say,
$$
f\in V_F,\,\,f'(z)\ne 0\,\forall\,z\in\mathbb{D}, {\rm meas}(\mathbb{D}-f(\mathbb{D}))=0\Rightarrow f\in {\rm Aut}(\mathbb{D}),
$$
where ${\rm meas}(A)$ is the Lebesgue measure of the Lebesgue measurable set $A$. Or maybe,
$$
f\in V_F,\,\,f'(z)\ne 0\,\forall\,z\in\mathbb{D}, \lim_{r\rightarrow 1^-}|f(re^{i\theta})|=1\,\,{\rm a.e.}\,\,{\rm in}\,\,\theta
\Rightarrow f\in {\rm Aut}(\mathbb{D}).
$$
This last question could be rephrased as follows:
$$
f\in V_F,\,\,f'(z)\ne 0\,\forall\,z\in\mathbb{D},\,f\,\,{\rm is}\,\,{\rm an}\,\,{\rm inner}\,\,{\rm function}
\Rightarrow f\in {\rm Aut}(\mathbb{D}).
$$

\section{The main result}

We can answer the two last questions that were raised in the previous section. We start by answering affirmatively the
last question.

\begin{theorem}
If $f\in V_F$, $f'(z)\ne 0\,\,\forall\,z\in\mathbb{D}$, $f$ is an inner function, then $f\in {\rm Aut}(\mathbb{D})$.
\end{theorem}
\noindent
{\bf Proof.} \\
We recall the following result, \\
\\
{\bf Theorem.} Every inner function is a uniform limit of Blaschke products. \\
\\
We refer to the theorem on page 175 of \cite{h}. Let $\{ B_n\}_{n=1}^{\infty}$ be a sequence of Blaschke products that uniformly
converge to $f$. Since $f\in V_F$ there exist a natural number $N_f$ such that $d_f(w)\le N_f\,\,\forall\,w\in\mathbb{D}$. By
Hurwitz Theorem we have $\lim_{n\rightarrow\infty} d_{B_n}(w)=d_f(w)\,\,\forall\,w\in\mathbb{D}$ and hence $d_{B_n}(w)=d_f(w)$
for $n\ge n_w$. We should note that $n_w$ depends on $w$ but it is a constant in a neighborhood of the point $w$. Hence the
Blaschke products in the tail subsequence $\{B_n\}_{n\ge n_w}$ all have a finite valence which is bounded from above by
$d_f(w)$ at the point $w$. In particular, the valence of these finite Blaschke products are bounded from above by the number
$N_f$ in definition 2.1 (of the set $V_F$). Hence we can extract a subsequence of these Blaschke products that have one and the 
same number of zeroes. Again by the Hurwitz Theorem it follows that $d_f(w)=N$ is a constant, independent of $w$, and so by
Theorem A in section 1 we conclude that $f(z)$ is a finite Blaschke product with exactly $N$ zeroes. By Theorem C in section 1
we conclude (using the assumption that $f'(z)\ne 0$, $\forall\,z\in\mathbb{D}$) that $f(z)\in {\rm Aut}(\mathbb{D})$. $\qed $ \\
\\
Next, we answer negatively the one before the last question.

\begin{theorem}
There exist functions $f\in V_F$ that satisfy $f'(z)\ne 0$ $\forall\,z\in \mathbb{D}$ and also ${\rm meas}(\mathbb{D}-f(\mathbb{D}))=0$
such that $f\not\in {\rm Aut}(\mathbb{D})$. In fact, we can construct such functions that will not be surjective and not injective.
\end{theorem}
\noindent
{\bf Proof.} \\
Consider the domain $\Omega=\mathbb{D}-\{x\in\mathbb{R}\,|\,0\le x<1\}$. Then $\Omega$ is the unit disc with a slit
along the non-negative $x$-axis. It is a simply connected domain. Let $g:\,\mathbb{D}\rightarrow\Omega$ be a Riemann mapping
(i.e. it is holomorphic and conformal. Finally, let $k\ge 2$ any natural integer and define $f=g^k$. This gives
the desired function. $\qed $ \\
\\
Can the result in Theorem 3.1 be generalized to higher complex dimensions? We make the obvious:

\begin{definition}
Let $V_F(n)$ ($n\in \mathbb{Z}^+$) be the set of all the holomorphic $f:\,\mathbb{D}^n\rightarrow\mathbb{D}^n$, such
that $\exists\,N=N_f$ (depending on $f$) for which $d_f(w)=|\{z\in\mathbb{D}^n\,|\,f(z)=w\}|$, $w\in\mathbb{D}^n$ satisfies
$d_f(w)\le N_f$ $\forall\,a\in\mathbb{D}^n$.
\end{definition}
\noindent
We ask if the following assertion holds true: \\
\\
If $f\in V_F(n)$ satisfies $\det J_f(z)\ne 0$ $\forall\,z\in\mathbb{D}^n$ and also $\lim_{r\rightarrow 1^-}|f(re^{i\theta_1},\ldots,re^{i\theta_n})|=1$
a.e. in $(\theta_1,\ldots,\theta_n)$ then $f\in {\rm Aut}(\mathbb{D}^n)$.

\noindent
{\it Ronen Peretz \\
Department of Mathematics \\ Ben Gurion University of the Negev \\
Beer-Sheva , 84105 \\ Israel \\ E-mail: ronenp@math.bgu.ac.il} \\ 
 
\end{document}